\documentclass[11pt,a4paper]{article}
\usepackage{amsfonts, amsmath, amssymb,  amsthm, color}
\usepackage{latexsym}
\title{\vspace{-2cm}
Scale analysis for an atmosphere flow
}
\author{\textbf{Donatella Donatelli}\\
        {\small Department of Information Engineering, Computer Science and Mathematics}\\
       {\small University of L'Aquila}\\
       {\small 67100 L'Aquila, Italy}\\
        $\scriptstyle\mathtt{donatella.donatelli@univaq.it}$
}

         \date{}
\newcommand{\e}{\varepsilon}		       
\newcommand{\R}{\mathbb{R}}

\newcommand{\ue}{u^{\varepsilon}}
\newcommand{\bu}{{\bf u}}
\newcommand{\bue}{{\bu^{\varepsilon}}}
\newcommand{\bued}{{\bu^{\varepsilon,\delta}}}

\newcommand{\bZ}{{\bf Z}}
\newcommand{\bZed}{{\bZ^{\varepsilon,\delta}}}

\newcommand{\pe}{p^{\varepsilon}}
\newcommand{\ped}{p^{\varepsilon, \delta}}
\newcommand{\omed}{\omega^{\varepsilon, \delta}}

\newcommand{\dive}{\mathop{\mathrm {div}}}
\newcommand{\curl}{\mathop{\mathrm {curl}}}

\newtheorem{theorem}{Theorem}[section]   

\newtheorem{lemma}[theorem]{Lemma}

\theoremstyle{definition}
\newtheorem{definition}[theorem]{Definition}
\newtheorem{remark}[theorem]{Remark}

\begin{document}

\maketitle
\begin{abstract}
In this paper we consider  a fluid dynamic model which describes  the atmospheric flow and we perform an asymptotic analysis for different time and length scales.
In particular we will focus on the two following cases:  when the Mach number dominates on the Rossby number yielding an incompressible regime or  when the  high  rotating (low Rossby number) and incompressibility regime act on the same scale ending up with a geostrophic balance. The limit analysis is performed in the ill prepared data setting and therefore we  will develop a rigorous and detailed analysis of the local decay and dispersive behavior of the related acoustic  waves. Finally we point out, that the set of equations analyzed in the paper may also fit in the artificial compressibility approximation methods.

\medbreak 
\medbreak
\end{abstract}
\tableofcontents
\newpage
\section{Introduction}
In this paper we perform a rigorous analysis of the singular limit as $\e\to 0$ of the following fluid dynamic  model,
\begin{equation}
\begin{cases}
\displaystyle{\partial_{t}\bue+\frac{1}{\e}({\bf g}\times \bue)+\frac{1}{\e^{2\beta}}\nabla \pe=\mu\Delta \bue-\left(\bue\cdot\nabla\right)\bue -\frac{1}{2}(\dive \bue)\bue}\\ \\
\e^{2\beta}\partial_{t}\pe+ \dive \bue=0,
\end{cases}
\label{i1}
\end{equation}
where  $x\in \Omega\subset \R^{3}$, $t\geq 0$, ${\bf u}$ is the fluid velocity and $p$ is the pressure, ${\bf g}=(0,0,1)$ is the rotation axis parallel to the the vertical variable $x_{3}$ and $\beta\geq 1$ or $\beta=1/2$. In the next section we will give  more details on  the geometry of the physical space $\Omega$.
\subsection{Motivation}
\label{Motivation}
The interest in studying the system \eqref{i1} comes from the fact that, according to the physical model we are considering,  it can be regarded in various ways. 
We mention here two examples to clarify this issue.

 As is well known, a rotating fluid is characterized by the Rossby number $\mathcal{R}$ that takes into account the Coriolis effect and, for example, it can be  used to model the movements of the oceans or of the atmosphere. If we denote by $U$ and $L$ the characteristic velocity and length scale of the fluid and by $f$ the local vertical component of the earth's rotation we have that the Rossby number is given by the formula
$$\mathcal{R}=\frac{U}{2fL}.$$ 
Performing the limit as $\mathcal{R}\to 0$ means that the scale motion of the fluid is much smaller than that of the earth. A very simple model for a rotating fluid is given by the incompressible Navier Stokes equations where  the external forces  are given by the rotating terms. Even if this type of model is quite simple from a mathematical point of view (we know a lot of results concerning the existence of solutions, see \cite{LPL96}), the incompressibility constraint from a computational point of view is very expensive. In fact discretization errors accumulate at each iteration and after a significant amount of error accumulation, the approximating algorithm breaks down. Moreover,  the incompressibility assumption may be not very realistic in modeling many physical phenomena. A way to overcome this troubles and, at the same time, to have a simple model, is to consider the compressible Navier Stokes system with density $\varrho$ and the pressure $p=p(\varrho)$ and to linearize it around a density constant state that for simplicity we can take as $\varrho=1$. Then, the linearized continuity equation assumes the form
\begin{equation}
\mathcal{M}\partial_{t} p+\dive\bu=0,
\label{mach}
\end{equation}
where $\mathcal{M}$ is the Mach Number given by
$$\mathcal{M}=\frac{U}{\sqrt{p'(1)}}.$$
The equation \eqref{mach} is a sort of ÒlinearizedÓ compressibility condition, namely
we  added to the incompressible constraint an ``artificial compressibility''  term $\mathcal{M}\partial_{t} p$ that vanishes as $\mathcal{M}\to 0$. This approximation method goes under the name of {\em artificial compressibility method} and was introduced by Chorin \cite{Ch68}, \cite{Ch69}, Temam \cite{Tem69a}, \cite{Tem69b} and Oskolkov \cite{Osk}.
If we take the   Rossby number  as $\mathcal{R}=\e$ and the  Mach number  as $\mathcal{M}=\e^{2\beta}$ then, the  above system is nothing else than the artificial compressibility approximation for a rotating fluid. It is important to remark that the first equation of the system \eqref{i1}  compared to  the balance momentum equation in the incompressible Navier Stokes equations has the extra term $-1/2(\dive\ue)\ue$ which has been added as a correction  to avoid the paradox of an increasing kinetic energy along the motion. 
Finally, \eqref{i1} can be considered as a family of perturbed systems, depending on a positive parameter $\e$, which may approximate in the limit the Navier Stokes equation.
In fact one of the main issue is to investigate in a rigorous way the limit as $\e\to 0$. The convergence of the  artificial compressibility system to the incompressible  Navier Stokes equation has been proved by  Temam \cite{Tem69a}, \cite{Tem69b} and his book \cite{Tem01}  on bounded domains, while in \cite{DM06}, \cite{DS11} and   \cite{DM10} in the case of the of  the whole space $\R^{3}$ and of the exterior domain and  then modified in a suitable way in \cite{Don08}  for the Navier Stokes Fourier system in $\R^{3}$, for the MHD system in \cite{Don13} and for the Navier-Stokes-Maxwell-Stefan system in \cite{DD18}. Here we have an additional difficulty in the limit analysis because of the rotating term,  we  will be more clear about this issue in the next Section \ref{difficulties}.

However, we can look at the system \eqref{i1} from a different point of view. As is well known, fluid dynamic equations are used to model various phenomena arising from physics, engineering, astrophysics. The model for atmosphere flows is given by the classical compressible fluid equation which include terms that take into account gravitation and rotation (see  \cite{K10}). One feature of the atmospheric flows is that they take place at different time and length scales and it is important to understand which phenomena occur according to the use of single scales or to the interactions of them (i.e. internal gravity waves, Rossby waves, cloud formation).
From a mathematical point of view, these various physical behavior give rise to different singular limits and to different asymptotic behaviors of the governing equations. Therefore  an other way to look at the  system \eqref{i1} is  to consider it as a simplified model that takes into account these different time and length scales for atmospheric flows. Indeed following \cite{K10} we set 
\begin{description}
\item[${\mathbf T=t_{ref}/\epsilon^{\alpha_{t}}}$] the characteristic time, with $\alpha_{t}>0$,
\item[${\mathbf L=l_{ref}/\epsilon^{\alpha_{x}}}$] the characteristic length $\alpha_{x}>0$,
\end{description}
where $\epsilon$ is the ratio between the thermal wind velocity $u_{ref}$ and the internal wave speed $c_{int}$, $l_{ref}=h_{sc}$, $t_{ref}=h_{sc}/u_{ref}$ and $h_{sc}$ is the density scale height. Then, a simplified set of equations (in the sense that we don't take into account temperature effects) that describes  the atmosphere flow is given by
\begin{equation*}
\frac{\epsilon^{\alpha_{t}}}{\epsilon^{\alpha_{x}}}\partial_{t}\bu+\bu\cdot\nabla\bu+\frac{{\bf g}\times \bu}{\epsilon^{\alpha_{x}-1}}+\frac{\epsilon^{\alpha_{\pi}}}{\epsilon^{3}}\nabla\tilde{\pi}=\bf{Q_{u}}
\end{equation*}
\begin{equation*}
\frac{\epsilon^{\alpha_{t}}}{\epsilon^{\alpha_{x}}}\partial_{t}\tilde{\pi}+\bu\cdot\nabla\tilde{\pi}+\frac{\gamma\pi}{\epsilon^{\alpha_{\pi}}}(\dive \bu)=\bf{Q_{\pi}},
\end{equation*}
where 
\begin{equation*}
\pi(t,x)=\bar{\pi}(x)+\epsilon^{\alpha_{\pi}}\Gamma \tilde{\pi}(t,x), \qquad \alpha_{\pi}>0.
\end{equation*}
and $\bf{Q_{u}}$ and $\bf{Q_{\pi}}$ are source terms and we assume that they have an appropriate high order in $\epsilon$ to not affect the leading order asymptotic analysis.

From the previous equations it is clear that according to the different values of $\alpha_{x}$, $\alpha_{t}$, $\alpha_{\pi}$ we obtain a hierarchy of models that corresponds to the undergoing physical process of multiple-scale regime. In particular, if we focus on the so called advection timescales, that is $\alpha_{x}=\alpha_{t}$, we can observe that for  $0\leq \alpha_{\pi}\leq 3$, $\alpha_{\pi}\neq 2$ we get as the leading order equations  the classical incompressible Navier Stokes equations.
On the other hand if, in  the advection time scale regimes we take $\alpha_{x}=\alpha_{t}=2$ and we balance the Coriolis and the pressure gradient terms by taking $\alpha_{\pi}=2$ we find out that the leading order dynamics is given by the geostrophic balance which corresponds to the quasi-geostrophic model in meteorology.
So we can conclude that our model \eqref{i1} is a sort of toy model for the previous cases, in particular the case $\beta\geq1$corresponds to  the first regime while $\beta=1/2$ is the geostrophic balance regime.

In the next section we will perform the formal limit for \eqref{i1}  in order to clarify  the different asymptotic behaviors that are a consequence of the different scaling regimes.

\subsection{Formal limit analysis and main mathematical difficulties}
\label{difficulties}
As already mentioned in the Section \ref{Motivation} in the system \eqref{i1} we can consider $\e^{2\beta}$ as the Mach Number $\mathcal{M}$ and $\e$ as the Rossby number $\mathcal{R}$ and we perform the limit as $\e\to 0$. Clearly, in this scenario, we have the competition of two effects that act simultaneously. If we consider the low Mach number limit which corresponds to the physical state in which the fluid speed is much smaller than the sound speed, the fluid density becomes constant, the velocity is soleinoidal and the fluid is incompressible. Low Rossby number corresponds to fast rotation and, from experimental data, at high rotating fluid becomes planar. Therefore we have to distinguish two cases according to the different values of $\beta$: 
\begin{itemize}
\item[\bf{[\boldmath$\beta\geq1$]}] As $\e\to 0$ first, the low Mach number regime dominates and the fluid becomes incompressible then, it stabilize to a planar flow and so we and up with an incompressible planar fluid which, if we denote by $\bu$ the limiting velocity, is described formally  by the set of equations
\begin{equation}
\partial_{t}\bu+\bu\cdot\nabla \bu-\Delta \bu=\nabla \pi, \quad \dive \bu=0.
\label{limite1}
\end{equation}
\item[\bf{[\boldmath$\beta=\frac{1}{2}$]}] As $\e\to 0$ the speed of rotation and the incompressibility act on the same scale, so the fluid becomes solenoidal and planar at the same time and we end up, at least formally, with a single linear equation given by
\begin{equation}{\bf g}\times\bu+\nabla \pi=0,\quad  \dive \bu=0.
\label{limite2}
\end{equation}
The equations \eqref{limite2} describe the geostrophic balance, but from them it is not possible to determine the time evolution of the flow, so we have to find a more complete description of the limiting behavior.
Formally, this can be done  by setting up an asymptotic expansion  and, by looking at high  order of $\e$, one obtains the following flow evolution  (for more details on this formal derivation see \cite{K10}),
$$
\partial_{t}(\nabla_{h}\pi-\pi)+\nabla_{h}^{\bot}\pi\cdot\nabla_{h}(\Delta_{h}\pi)=\Delta_{h}^{2}\pi.
$$

\end{itemize}

When we try to make rigorous the previous analysis, one of the main problems is that as $\e\to 0$ the velocity fields develops very fast oscillating waves in time (the so called {\em acoustic waves}). These waves are supported by the gradient part of the velocity field, they propagate along the motions and  give rise to the lost of compactness for the nonlinear terms. It is clear that if we take the initial data constructed in such a way that they are supported in the kernel of the acoustic wave operator, namely ``well-prepared initial data'' the  limiting process is not affected by these waves. 

In this paper we take very general initial data,  we only require the boundedness of the initial energy, hence, as a consequence we have to deal with the presence  of the acoustic waves. Given the particular geometry of the domain, we will work on an infinite slab,  we expect that at a certain time these waves will loose and disperse they energy in the space domain. So, in order to control the oscillations and gain some sort of compactness we will develop a rigorous and detailed analysis of the local decay and dispersive behavior of these waves. 

This type of analysis will be different according to the different values of $\beta$.  In fact we have the dispersive behavior of the acoustic waves and at the same time  the fluid is under the effects of the centrifugal force that becomes large as $\e\to 0$. For these reasons we cannot use the classical dispersive estimate of Strichartz type as in \cite{DM06}. 
To  be precise, as $\beta \geq 1$ the decay of the acoustic wave is strong enough to eliminate the centrifugal force and we will perform some rigorous decay estimate of the acoustic waves in the spirit of D'Ancona and Racke  \cite{DR12} and Sogge \cite{SmSo00}, see also \cite{DFN12}.  If $\beta=1/2$ then, the incompressible regime and the high rotation occur at the same scale, we cannot exploit the local decay of the oscillating waves but we have to analyze the spectral properties of the rotating operator. We have to  show that the fast oscillating parts of the gradient live in the space orthogonal to the kernel of the rotating operator and so they don't affect our limiting process. The basic tool in this case will be  the RAGE theorem, see  \cite{DFN12}

\subsection{Plan of the paper}
In the next Section \ref{Setting} we set the problem, we define the notion of weak solutions we are going to use and we state the main results of this paper, Theorem \ref{T1} and Theorem \ref{T2}. In Section \ref{SE} we perform the uniform a priori estimate (with respect to $\e$) satisfied by the solutions of our system \eqref{i1} for any $\beta\geq 1$ and $\beta=1/2$. In Section \ref{beta maggiore di uno} we perform the limit analysis for the case $\beta\geq 1$  and we prove  the Theorem \ref{T1}. Finally in Section \ref{beta uguale1/2} we deal with the case $\beta=1/2$ and we prove the Theorem \ref{T2}.

\section{Setting of the problem and main results}
\label{Setting}
\subsection{Notation}
Before setting up or problem we fix here the main notations we are going to use through the paper.
\begin{itemize}
\item $C^{\infty}_{0}([0,T)\times \Omega)$ is the space of $C^{\infty}$ functions with compact support
\item $W^{k,p}(\Omega)$ is the usual Sobolev space on $\Omega$ and $H^{k}(\Omega)=W^{k,2}(\Omega)$. 
\item The notations
 $L^{p}_{t}L^{q}_{x}$ and $L^{p}_{t}W^{k,q}_{x}$ will abbreviate respectively  the spaces $L^{p}([0,T);L^{q}(\Omega))$, and $L^{p}([0,T);W^{k,q}(\Omega))$.
\item $Q$ and $P$ are the  Leray's projectors  on the space of gradients vector fields and  on the space of divergence - free vector fields respectively, namely
$$ P=I-Q.$$
\item $\nabla^{\bot}_{h}f$ denotes the vector $(\partial_{x_{2}}f, -\partial_{x_{1}}f)$. 
The differential operators $\nabla_{h}$, ${\dive}_{h}$, $\Delta_{h}$ 
denote the usual $\nabla$, ${\dive}$, $\Delta$ applied on the horizontal variables $x_{h}=(x_{1}, x_{2})$.
\item For a function $f$ the vertical average on the one dimensional torus  ${\mathbb{T}}^{1}$ is defined  by 
$$\langle f(x_{h})\rangle=\frac{1}{|{\mathbb{T}}^{1}|}\int_{{\mathbb{T}}^{1}}f(x_{h},x_{3})dx_{3}.$$
\end{itemize}
\subsection{Setting of the problem}
We consider the following system
\begin{equation}
\begin{cases}
\displaystyle{\partial_{t}\bue+\frac{1}{\e}({\bf g}\times \bue)+\frac{1}{\e^{2\beta}}\nabla \pe=\mu\Delta \bue-\left(\bue\cdot\nabla\right)\bue -\frac{1}{2}(\dive \bue)\bue}\\ \\
\e^{2\beta}\partial_{t}\pe+ \dive \bue=0,
\end{cases}
\label{1}
\end{equation}
where $x\in \Omega\subset \R^{3}$, $t\geq 0$, $\mu\in \R$ and ${\bf g}=(0,0,1)$ is the rotation axis parallel to the the vertical variable $x_{3}$. The geometry of the physical space $\Omega$ is given by an infinite slab,
\begin{equation*}
\Omega=\R^{2}\times (0,1).
\label{2}
\end{equation*}
and we will denote the horizontal variable as $x_{h}=(x_{1}, x_{2})$. 

For the velocity field $\bue$ we assume the complete slip boundary conditions
\begin{equation}
\bue\cdot{\bf n}|_{\partial\Omega}=0,\qquad [\mathbb{S}{\bf n}]\times {\bf n}|_{\partial\Omega}=0,
\label{3}
\end{equation}
where  ${\bf n}$ denotes the outer normal vector to the boundary and $\mathbb{S}$ is the viscous stress tensor given by
$$\mathbb{S}(\nabla\bue)=\mu\left(\nabla\bue+\nabla^{t}\bue-\frac{2}{3}\dive\mathbb{I}\right), \quad \mu >0.$$ 
From now on, without loss of generality, for simplicity we set $\mu=1$.

In order do deal with the boundary conditions \eqref{3} it is more convenient to reformulate the state variable in a periodic setting for the variable $x_{3}$. In fact we will take
\begin{equation*}
\Omega=\R^{2}\times \mathbb{T}^{1},
\label{4}
\end{equation*}
where $\mathbb{T}^{1}$ is the one dimensional torus and where the pressure is extended even in the third variable,
$$\pe(x_{1}, x_{2}, -x_{3})=\pe(x_{1}, x_{2}, x_{3}),$$
as well as the horizontal component of the velocity ${\bu}^{\e}_{h}=(\ue_{1},\ue_{2})$,
$$\ue_{j}(x_{1}, x_{2}, -x_{3})=\ue_{j}(x_{1}, x_{2}, x_{3}),\quad j=1,2,$$
while the vertical component $\ue_{3}$ is taken odd,
$$\ue_{3}(x_{1}, x_{2}, -x_{3})=-\ue_{3}(x_{1}, x_{2}, x_{3}).$$
Furthermore we assign to the system \eqref{1} the following  initial conditions
\begin{equation}
\bue(x,0)=\bu^{\e}_{0}(x),\  \pe(x,0)=\pe_{0}(x).\\
 \label{initial data}
\end{equation}
The  regularity and the limiting behavior as $\e\to 0$ of the initial data  \eqref{initial data} deserve a little discussion. Indeed  the system \eqref{1} requires the initial conditions \eqref{initial data} while the target equations \eqref{limite1} and \eqref{limite2} require only the  initial condition for the velocity $\bu$. Hence, our approximation will be consistent if the initial datum on the pressure $\pe$  will be eliminated by an ``initial layer'' phenomenon. Since in the limit we have to deal with weak solutions in the sense of the Definition \ref{weaksol} it is reasonable to require the finite energy constraint to be satisfied by the approximating sequences $(\bue, \pe)$. So we can deduce a natural behavior to be imposed on the initial data $(\bu^{\e}_{0},  \pe_{0})$, namely
\begin{equation}
\begin{split}
\bu^{\e}_{0}, \quad& \pe_{0} \in L^{2}(\Omega)\\
\bu^{\e}_{0}\rightharpoonup \bu_{0}, \quad&  \pe_{0}\rightharpoonup p_{0}\quad \text{weakly in}\  L^{2}(\Omega).
\end{split}
\label{ID}
\end{equation}

For completeness we recall the notion of weak solutions for the system \eqref{1} we are going to use.
\begin{definition}
\label{weaksol}
We say that  a pair $\bue, \pe$ is a  weak solution to the system \eqref{1} in $(0,T)\times \Omega$ if $\ue\in L^{\infty}([0,T];L^{2}(\R^{3}))\cap L^{2}([0,T];\dot H^{1}(\R^{3})) $.
$\pe \in L^{\infty}([0,T];L^{2}(\R^{3}))$ and   they 
satisfy  \eqref{1} in the sense of distributions, namely
\begin{align}
\int_{0}^{T}\!\!\int_{\Omega}\Big(\bue
\partial_{t}\varphi &-((\bue\cdot\nabla) \bue+\frac{1}{2}(\dive\bue)\bue)\cdot\varphi\notag\\&-\frac{1}{\e}({\bf g}\times \bue)\cdot\varphi+\frac{1}{\e^{2\beta}}\pe\dive\varphi \Big) dxdt \notag\\
&=\int_{0}^{T}\!\!\int_{\Omega}\nabla\bue:\nabla\varphi dxdt-\int_{\Omega}\bu^{\e}_{0}\cdot \varphi(0,\cdot) dx,
\label{weak}
\end{align}
for any $\varphi\in C^{\infty}_{0}([0,T)\times \Omega;\R^{3})$ and
\begin{align*}
&\int_{0}^{T}\!\!\int_{\Omega}\left(\e^{2\beta}\pe\partial_{t}\varphi +\bue\cdot\nabla\varphi\right )dxdt=
-\int_{\Omega}\e^{2\beta}\pe_{0}\varphi (0,\cdot)dx,
\end{align*}
for any $\varphi\in C^{\infty}_{0}([0,T)\times \Omega)$. Moreover the following energy inequality holds
\begin{align*}
\frac{1}{2}\int_{\Omega}&(|\bue(x,t)|^{2}+|\pe(x,t)|^{2})dx+\int_{0}^{t}\!\!\int_{\Omega}|\nabla \bu(x,s)|^{2}dxds\\ \leq
&\frac{1}{2}\int_{\Omega}(|\bu^{\e}_{0}(x)|^{2}+|\pe_{0}(x)|^{2})dx ,\qquad \text{for all $t\geq 0$}.
\end{align*}
\end{definition}
The proof of the existence of global in time weak solutions  for \eqref{1}  is omitted since, in the spirit of   Temam (see Chapter III, Theorem 8.1 in \cite{Tem01}),   it follows by standard finite dimensional Galerkin type approximations with the necessary modification due to the domain $\Omega$.

\subsection{Main Results}
Now we are ready to  state the  main results of this paper.

\begin{theorem}[{\bf Case \boldmath$\beta\geq 1$}]
\label{T1}
Assume that $\bue, \pe$ are weak solutions of the system \eqref{1} with initial data \eqref{initial data} satisfying \eqref{ID}, then there exists $\bu\in L^{2}(0,T;W^{1,2}(\Omega))$,  
such that
\begin{equation}
\bue\rightharpoonup \bu \quad \text{weakly in  $L^{2}(0,T;W^{1,2}(\Omega)),$}
\label{r1bis}
\end{equation}
\begin{equation}
\bue\longrightarrow \bu \quad \text{strongly in  $L^{2}_{loc}((0,T)\times \Omega)$},
\label{r2bis}
\end{equation}
where $\bu=[\bu_{h}(t,x_{h}),0]$  is the unique weak solution of the 2D incompressible Navier Stokes equation 
\begin{equation}
{\dive}_{h} \bu=0,
\label{r4bis}
\end{equation}
\begin{equation}
\partial_{t}\bu+(\bu\cdot\nabla_{h})\bu+\nabla_{h} \pi=\Delta_{h} \bu.
\label{r5}
\end{equation}
\end{theorem}

\begin{theorem}[{\bf Case \boldmath$\beta=1/2$}]
\label{T2}
Assume that $\bue, \pe$ are weak solutions of the system \eqref{1} with initial data \eqref{initial data} satisfying \eqref{ID}, then there exists $\bu\in L^{2}(0,T;W^{1,2}(\Omega))$,  $\pi\in L^{\infty}(0,T;L^{2}(\Omega))$such that
\begin{equation}
\bue\rightharpoonup \bu \quad \text{weakly in  $L^{2}(0,T;W^{1,2}(\Omega)),$}
\label{r1}
\end{equation}
\begin{equation}
\bue\longrightarrow \bu \quad \text{strongly in  $L^{2}_{loc}((0,T)\times \Omega)$},
\label{r2}
\end{equation}
\begin{equation}
\pe\rightharpoonup \pi \quad \text{weakly in  $L^{\infty}(0,T;L^{2}(\Omega))$},
\label{r3}
\end{equation}
where $\bu$ and $\pi$ satisfy
\begin{equation}
{\dive}_{h} \bu=0,
\label{r4}
\end{equation}
\begin{equation}
{\bf g} \times \bu+\nabla \pi=0
\label{r5}
\end{equation}
and $\pi$ is a solution in the sense of distribution of the equation
\begin{equation}
\label{r6}
\partial_{t}(\nabla_{h}\pi-\pi)+\nabla_{h}^{\bot}\pi\cdot\nabla_{h}(\Delta_{h}\pi)=\Delta_{h}^{2}\pi.
\end{equation}
\end{theorem}
Concerning the limiting equations \eqref{r5} and \eqref{r6} it is important to remark that in the geostrophic balance regime the acoustic waves are asymptotically filtered out. This type of phenomenon corresponds to quasi-geostrophic regime with the advection timescale.


\section{Energy estimate and uniform bounds}
\label{SE}

We define the energy functional associated to the system \eqref{1} as 
\begin{equation*}
\label{E}
E(t)=\frac{1}{2}\int_{\Omega}\left( |\bue(x,t)|^{2}+ |\pe(x,t)|^{2}\right)dx.
\end{equation*}
By standard computations it is straightforward to prove that  the weak solutions of the system \eqref{1}, satisfy  the energy equality

\begin{equation}
\label{E2}
E(t)+\int_{0}^{t}\!\!\int_{\R^{3}}|\nabla \bue(x,s)|^{2}dxds\leq E(0).
\end{equation}
As a consequence  \eqref{E2} we obtain the following uniform bounds

\begin{equation}
\bue,\ \pe \quad  \text{are bounded in} \ L^{\infty}([0,T];L^{2}(\Omega)),\label{5}
 \end{equation}
  \begin{equation}
 \nabla\bue \quad  \text{is bounded in $L^{2}([0,T]\times\Omega).$}  \label{6}
 \end{equation}
 By combining \eqref{E2} with  standard Sobolev embeddings  we deduce that
  \begin{equation}
 \bue \quad  \text{is bounded in $L^{\infty}([0,T];L^{2}(\Omega))\cap L^{2}([0,T];L^{6}(\Omega)).$}  \label{7}
 \end{equation}
Using together \eqref{6} and \eqref{7} we have that 
  \begin{equation}
  \begin{split}
&(\bue \cdot \nabla)\bue\quad  \bue\dive\bue\\
 &\text{are bounded in $L^{2}([0,T];L^{1}(\Omega))\cap L^{1}([0,T];L^{3/2}(\Omega)).$}  
\end{split}
\label{8}
\end{equation}

We point out that the previous estimates are uniform in $\e$ and hold for any value of $\beta>0$.

\section{Case $\beta\geq 1$}
\label{beta maggiore di uno}
This section is devoted to the analysis of the limiting behavior as $\e\to 0$ in the case $\beta\geq 1$.  As a first step we  recover  the convergence results that follows from the bounds of the previous Section \ref{SE}. Then, as already mentioned in the introduction we have to deal with the  high oscillations of the acoustic waves. Since these waves are supported by the gradient part of the velocity, we decompose the velocity in its gradient and soleinoidal part and, since the  domain $\Omega$ is infinite, we will be able to get decay estimates for the acoustic potential (gradient part of the velocity). As a last step we deal with the convergence of the velocity divergence  free part and, in order to get the strong convergence, we have to estimate  the vertical average and the oscillations of the solenoidal component of $\bue$.

\subsection{First convergence results}
\label{first conv}
From \eqref{6} and \eqref{7} we have that 
\begin{equation}
\label{9}
\bue \rightharpoonup \bu \quad \text{weakly in $L^{2}([0,T];H^{1}(\Omega)) $}.
\end{equation}
Hence, by taking into account \eqref{5} and \eqref{9}, and letting $\e\to 0$ in $\eqref{1}_{2}$ we obtain
\begin{equation}
\label{10}
\dive \bu =0 \quad \text{a.e. in $(0,T)\times \Omega$}.
\end{equation}
Moreover, if we apply the Leray projector $P$ on $\eqref{1}_{1}$, as   $\e\to 0$ we have 
\begin{equation*}
\label{11}
P({\bf g}\times \bu)=0,
\end{equation*}
from which it follows that ${\bf g}\times \bu= \nabla G$, for a certain potential $G$. As a consequence, the limiting velocity horizontal component $\bu_{h}=(u_{1}, u_{2})$ does not depend on the vertical variable $x_{3}$, then  by using  \eqref{10},  we have that $\partial_{x_{3}}u_{3}=0$. If we take into account the boundary condition \eqref{3} and the fact that $\bu\in L^{2}_{t}L^{2}_{x}$ we can conclude that
\begin{equation}
\label{12}
\bu=(\bu_{h}(t, x_{h}),0), \quad u_{3}=0.
\end{equation}
Now, if we choose a test function of the form $\varphi(x,t)=(\varphi_{h}(t, x_{h}),0)$, $\dive \varphi=0$, it is possible to  pass into the limit in the weak formulation of $\eqref{1}_{1}$, provided we know how to handle  the nonlinear terms $(\bue\cdot\nabla) \bue$,  $\bue\dive \bue$. In the next section,  by studying separately the  gradient and solenoidal part of $\bue$, we will get stronger  convergence  results.


\subsection{Acoustic equation and estimates} 
In this section we will analyze the behavior of the acoustic waves in order to control their fast oscillation in time. To this end we rewrite the system \eqref{1} in the following form
\begin{equation}
\begin{cases}
\e^{2\beta}\partial_{t}\pe+ \dive \bue=0\\ \\
\displaystyle{\e^{2\beta}\partial_{t}\bue+ \nabla \pe=-\e^{2\beta -1}({\bf g}\times \bue)+\e^{2\beta}(\mu\Delta \bue-\left(\bue\cdot\nabla\right)\bue -\frac{1}{2}(\dive \bue)\bue)}.
\end{cases}
\label{ae}
\end{equation}
We can observe that the underlying structure of the system \eqref{ae} is that of a wave equation in fact  it goes under the name of  \emph{acoustic wave system}. In order to simplify the computations we rewrite \eqref{ae} as
\begin{equation}
\label{ae2}
\begin{cases}
\e^{2\beta}\partial_{t}\pe+ \dive \bue=0\\ \\
\e^{2\beta}\partial_{t}\bue+ \nabla \pe=-\e^{2\beta -1}\mathbb{F}^{\e}_{1}+\e^{2\beta}\dive \mathbb{F}^{\e}_{2}+\e^{2\beta}\mathbb{F}^{\e}_{3},
\end{cases}
\end{equation}
where by taking into account \eqref{6}-\eqref{8} we have $\mathbb{F}^{\e}_{1}\in L^{\infty}_{t}L^{2}_{x}$, $\mathbb{F}^{\e}_{2}\in L^{2}_{t}L^{2}_{x}$, $\mathbb{F}^{\e}_{3}\in L^{2}_{t}L^{1}_{x}$. 
Since the equations in \eqref{ae2} are satisfied only in a weak sense it is more convenient to regularize them in order to deal with smooth solutions.
To this purpose, given a function $v$, and $j_{\delta}$ a standard Friedrich's mollifier we denote by
$$v^{\delta}=j_{\delta}\ast v$$
the regularized  function. Having in mind these notations we regularize the system \eqref{ae2} and we get

\begin{equation}
\begin{cases}
\e^{2\beta}\partial_{t}\ped+ \dive \bued=0\\ \\
\e^{2\beta}\partial_{t}\bued+ \nabla \ped=\e^{2\beta -1}\mathbb{F}^{\e, \delta}_{1}+\e^{2\beta}\dive \mathbb{F}^{\e,\delta}_{2}+\e^{2\beta}\mathbb{F}^{\e,\delta}_{3},
\end{cases}
\label{ae33}
\end{equation}
where, for the right hand side of $\eqref{ae33}_{2}$ the following uniform  bounds in $\e$ hold
\begin{equation*}
\label{13}
\|\mathbb{F}^{\e, \delta}_{1}\|_{L^{2}_{t}H^{k}}+\|\mathbb{F}^{\e, \delta}_{2}\|_{L^{2}_{t}H^{k}}+\|\mathbb{F}^{\e, \delta}_{3}\|_{L^{2}_{t}H^{k}}\leq c(k,\delta),\quad \text{for any $k=0,1,\ldots$}
\end{equation*}
Since the acoustic waves are supported by the gradient part of the velocity fields,  we decompose $\bued$  in the following way
\begin{equation}
\label{14}
\bued=\bZed+\nabla \Psi^{\e,\delta},
\end{equation}
where 
$\bZed=P\bued$,  $\nabla \Psi^{\e,\delta}=Q\bued$ and we rewrite \eqref{ae33} in terms of the acoustic potential $\Psi^{\e,\delta}$,

\begin{equation}
\e^{2\beta}\partial_{t}\ped+\Delta \Psi^{\e,\delta}=0,
\label{ae3}
\end{equation}
\begin{equation}
\e^{2\beta}\partial_{t}\Psi^{\e,\delta}+\ped\!=\e^{2\beta -1}\Delta^{-1}\!\dive\mathbb{F}^{\e, \delta}_{1}+\e^{2\beta}\Delta^{-1}\!\!\dive(\dive \mathbb{F}^{\e,\delta}_{2}\!+\Delta^{-1}\!\dive\mathbb{F}^{\e,\delta}_{3}).
\label{ae4}
\end{equation}
It is clear now that  to estimate $\Psi^{\e,\delta}$ we have to exploit the dispersive behaviour of the system \eqref{ae3}, \eqref{ae4} from which we will   deduce the local decay of the acoustic potential. So we recall here the following lemma for the proof of which see Feireisl et al. \cite{FGVN12} or D'Ancona and Racke \cite{DR12}
\begin{lemma}
\label{lemmadisp}
Consider $\varphi\in C^{\infty}_{0}(\R^{2})$. Then we have
\begin{equation}
\label{l1}
\int_{-\infty}^{\infty}\int_{\Omega}\left|\varphi(x_{h})\exp\left({\rm i}\sqrt{-\Delta}t\right)[v]\right|^{2} dxdt\leq c(\varphi)\|v\|^{2}_{L^{2}(\Omega)}.
\end{equation}
Moreover, on any compact set $K\subset\Omega$, $m>0$, we have
\begin{align}
\int_{0}^{T}\!\!\int_{K}&\left|\exp\left({\rm i}\sqrt{-\Delta}\frac{t}{\e^{m}}\right)[v]\right|^{2}dxdt\notag\\
&\leq\e^{m}\int_{0}^{\infty}\int_{\Omega}\left|\exp\left({\rm i}\sqrt{-\Delta}t\right)[v]\right|^{2} dxdt\leq \e^{m}c\|v\|^{2}_{L^{2}(\Omega)},
\label{l2}
\end{align}
and 
\begin{align}
\int_{0}^{T}\!\!\int_{K}&\left|\int_{0}^{t}\exp\left({\rm i}\sqrt{-\Delta}\frac{t-s}{\e^{m}}\right)[g(s)]\right|^{2}dxdt\notag\\
&\leq cT\e^{m}\int_{0}^{T}\left\|\exp\left({\rm i}\sqrt{-\Delta}\frac{s}{\e^{m}}\right)[g(s)]\right\|^{2}_{L^{2}(\Omega)}=\e^{m}\|g\|_{L^{2}((0,T)\times \Omega)}.
\label{l3}
\end{align}
\end{lemma}
By means of Duhamel's formula the gradient of the solution $\Psi^{\e,\delta}$ of  \eqref{ae4} is given by 
\begin{align*}
\label{15}
 &\nabla\Psi^{\e,\delta}(t) \\
 &= \frac{1}{2}
\exp \left( {\rm i} \sqrt{ -\Delta} \frac{t}{\e^{2\beta}} \right) \left[ \nabla\Psi^{\e, \delta}(0) + \frac{{\rm i}}{\sqrt{ -\Delta}}\nabla\ped(0)\right]\\
&+\frac{1}{2}
\exp \left( -{\rm i} \sqrt{ -\Delta} \frac{t}{\e^{2\beta}} \right) \left[ \nabla\Psi^{\e,\delta}(0) - \frac{{\rm i}}{\sqrt{ -\Delta}}\nabla\ped(0)\right]\notag\\
& +\frac{\e^{-1}}{2}\int_{0}^{t}\!\left(\exp \left( {\rm i} \sqrt{ -\Delta} \frac{t-s}{\e^{2\beta}} \right)+\exp \left(- {\rm i} \sqrt{ -\Delta} \frac{t-s}{\e^{2\beta}} \right)\!\!\right)\![\nabla\Delta^{-1}\!\dive\mathbb{F}^{\e, \delta}_{1}]ds\notag\\
& +\frac{1}{2}\int_{0}^{t}\!\left(\exp \left( {\rm i} \sqrt{ -\Delta} \frac{t-s}{\e^{2\beta}} \right)+\exp \left(- {\rm i} \sqrt{ -\Delta} \frac{t-s}{\e^{2\beta}} \right)\!\!\right)\![\nabla\Delta^{-1}\!\dive\!\dive \mathbb{F}^{\e,\delta}_{2}]ds\notag\\
& +\frac{1}{2}\int_{0}^{t}\!\left(\exp \left( {\rm i} \sqrt{ -\Delta} \frac{t-s}{\e^{2\beta}} \right)+\exp \left(- {\rm i} \sqrt{ -\Delta} \frac{t-s}{\e^{2\beta}} \right)\!\!\right)\![\nabla\Delta^{-1}\!\dive\mathbb{F}^{\e,\delta}_{3}]ds.\notag
\end{align*}
Now, by applying \eqref{l2} and \eqref{l3} with $m=2\beta$ we have the following uniform in $\e$ decay estimate for the acoustic potential,
\begin{equation}
\label{16}
\int_{0}^{T}\!\!\!\|\nabla\Psi^{\e,\delta}\|^{2}_{L^{2}(K)}dt\leq (\e^{2\beta-1}\!+\e^{2\beta})c(\delta, K,T), 
\end{equation} 
for any compact set $K\subset\Omega$ and  $\beta\geq 1$. Hence, we can conclude that the effects of the acoustic potential vanishes as soon as $\e\to 0$, so the limiting behaviour of the system \eqref{1} depends only on the soleinoidal component of the velocity field. Finally, we point out that in order to have a negligible effect of the acoustic potential  it is enough to require $\beta >1/2$.

\subsection{Convergence of the soleinoidal part of the velocity}
From the previous section we understood that the asymptotic behavior of the system \eqref{1} depends on the solenoidal part of $\bued$, hence in  this section we will  inquire on $P\bued$.
Since in the Section \ref{first conv}  we have proved that $\bued$ converges to a function $\bu$ which depends only on the horizontal variables, in order  get the strong  convergence of  $P\bued$  a first step is to establish the compactness of the vertical average of $\bued$. Therefore we rewrite the equation for $\bued$ as follows
\begin{equation}
\label{17}
\e\partial_{t}\bued+({\bf g}\times \bued)=\e{\bf S}^{\e,\delta}-\e^{1-2\beta} \nabla \ped,
\end{equation}
where
$${\bf S}^{\e,\delta}=\mu\Delta \bued-(\bued\cdot\nabla)\bued -\frac{1}{2}(\dive \bued)\bued$$
and ${\bf S}^{\e,\delta}$ is bounded in $L^{2}_{t}H^{k}_{x}$, for any fixed $k$ and $\delta$.
We take now the vertical average of \eqref{17},
\begin{equation}
\label{18}
\e\partial_{t}\langle\bued\rangle+({\bf g}\times \langle\bued\rangle)=\e\langle {\bf S}^{\e,\delta}\rangle-\e^{1-2\beta} \nabla \langle\ped\rangle.
\end{equation}
Since $\bZed$ is soleinoidal   one can easily check that $P({\bf g}\times \langle\bZed\rangle)=0$. Using a test function $\varphi\in C^{\infty}_{0}(\Omega;\R^{3})$, $\dive \varphi=0$ in the weak formulation of \eqref{18} we obtain
\begin{equation*}
\label{19}
\partial_{t}\int_{\Omega}\langle\bued\rangle\cdot \varphi dx=\int_{\Omega}\langle {\bf S}^{\e,\delta}\rangle\cdot\varphi dx-\frac{1}{\e}\int_{\Omega}({\bf g}\times \langle\nabla\Psi^{\e,\delta}\rangle)\cdot\varphi dx.
\end{equation*}
If we use \eqref{16} and taking into account that $\beta\geq 1$, by applying Lions Aubin lemma arguments (see \cite{Si}) we can conclude that
\begin{equation}
\label{20}
\langle\bZed\rangle \longrightarrow \bu^{\delta}, \quad\text{strongly in $L^{2}((0,T)\times K)$},
 \end{equation}
for any compact set $K\subset \Omega$ and any fixed $\delta$. It is worthful to remark that at this step it is crucial the estimate \eqref{16} and here we need the restriction on the values of $\beta$ ($\beta\geq 1$). Finally, 
 in order to get the limiting behaviour of $\bZed$ it is fundamental to  control some possible  oscillations.
 
Since we proved that the horizontal component of $\bued$ is compact we can infer that the oscillations are due to the vector fields that depends on $x_{3}$. In the remaining part of this section we will study and estimate these oscillations and we will show that they don't interfere in the convergence of the nonlinear terms.  For any function $f$ we denote the oscillation as
$$\{f\}(x)=f(x)-\langle f\rangle(x_{h}).$$
Notice that $\{f\}(x)$ has zero vertical mean and so it can be written for some function $I$ as 
$$\{f\}(x)=\partial_{x_{3}}I(x),\quad \text{with} \quad \int_{\mathbb{T}^{1}}I(x)dx_{3}=0.$$
Then we define for any $i,j=1,2,3$
$$\omed_{i,j}=\partial_{x_{i}}Z^{\e,\delta}_{j}-\partial_{x_{j}}Z^{\e,\delta}_{i}=\partial_{x_{i}}u^{\e,\delta}_{j}-\partial_{x_{j}}u^{\e,\delta}_{i}.$$
From  \eqref{18} we have that $\omed_{i,j}$ satisfy the following equations 
\begin{equation}
\label{21}
\e\partial_{t}\omed_{1,2}+{\dive}_{h}[\bZed]_{h}=\e\left(\partial_{x_{1}}S^{\e,\delta}_{2}- \partial_{x_{2}}S^{\e,\delta}_{1}\right)-\Delta_{h}\Psi^{\e,\delta},
\end{equation}

\begin{equation}
\label{22}
\e\partial_{t}\omed_{1,3}+\partial_{x_{3}}Z^{\e,\delta}_{2}=\e\left(\partial_{x_{1}}S^{\e,\delta}_{3}- \partial_{x_{3}}S^{\e,\delta}_{1}\right)-\partial^{2}_{x_{3}x_{2}}\Psi^{\e,\delta},
\end{equation}

\begin{equation}
\label{23}
\e\partial_{t}\omed_{2,3}-\partial_{x_{3}}Z^{\e,\delta}_{1}=\e\left(\partial_{x_{2}}S^{\e,\delta}_{3}- \partial_{x_{3}}S^{\e,\delta}_{2}\right)-\partial^{2}_{x_{3}x_{1}}\Psi^{\e,\delta}.
\end{equation}

Now by using  the decomposition \eqref{14} we rewrite  the nonlinear terms of the system \eqref{1} as,  
\begin{align}
(\bued\cdot\nabla)&\bued+\frac{1}{2}\bued\dive\bued=\dive(\bued\otimes\bued)-\frac{1}{2}\bued\dive\bued\notag\\
&=\dive(\bZed\otimes \bZed)+\dive(\nabla\Psi^{\e,\delta}\otimes \nabla\Psi^{\e,\delta})\notag\\
&+\dive(\bZed\otimes \nabla\Psi^{\e,\delta})+\dive(\nabla\Psi^{\e,\delta}\otimes \bZed)-\frac{1}{2}\bued\Delta \Psi^{\e,\delta}.
\label{24}
\end{align}
By taking into account  the local decay \eqref{16} of the acoustic potential we see that the only term of \eqref{24} who requires a detailed analysis  is
\begin{equation}
\label{25}
\dive(\bZed\otimes \bZed)=\frac{1}{2}\nabla|\bZed|^{2}-\bZed\times\curl[\bZed].
\end{equation}
We focus on the second term of right-hand side of \eqref{25}, 
\begin{align}
\label{26}
\bZed\times\curl[\bZed]=&\langle\bZed\rangle\times\curl\langle\bZed\rangle\notag
\\&+\partial_{x_{3}}\left(\langle\bZed\rangle\times\curl I[\bZed]+I [\bZed]\times  \curl\langle\bZed\rangle\right)\notag\\
&+\partial_{x_{3}}I [\bZed]\times\partial_{x_{3}}\curl\langle\bZed\rangle.
\end{align}
The first term of the right hand side of \eqref{26} is compact because of \eqref{20}, the second term has zero vertical mean, hence we have to study carefully only the last one.  For any $j=1,2,3 $ 
we have
\begin{align}
\label{27a}
&[\partial_{x_{3}}I [\bZed]\times\partial_{x_{3}}\curl\langle\bZed\rangle]_{j}\notag\\
&=\partial_{x_{3}}I [Z^{\e,\delta}_{i}]\partial_{x_{3}}(\partial_{x_{i}}I [Z^{\e,\delta}_{j}]-\partial_{x_{j}}I [Z^{\e,\delta}_{i}]=\partial_{x_{3}}I [Z^{\e,\delta}_{i}]\partial_{x_{3}}I [\omega^{\e,\delta}_{i,j}].
\end{align}
From the relations  \eqref{21} -\eqref{23} it is easy to obtain, 
\begin{equation}
\label{22b}
\e\partial_{t}(\partial_{x_{3}}I[\omed_{1,3}])+\partial_{x_{3}}^{2}I[Z^{\e,\delta}_{2}]\!=\!\e\!\left(\partial_{x_{1}}(S^{\e,\delta}_{3}-\!\langle S^{\e,\delta}_{3}\rangle)\!-\! \partial_{x_{3}}S^{\e,\delta}_{1}\right)-\partial^{2}_{x_{3}x_{2}}\!\Psi^{\e,\delta},
\end{equation}

\begin{equation}
\label{23b}
\e\partial_{t}(\partial_{x_{3}}I[\omed_{2,3}])- \partial_{x_{3}}^{2}I[Z^{\e,\delta}_{1}]\!=\!\e\!\left(\partial_{x_{2}}(S^{\e,\delta}_{3}-\!\langle S^{\e,\delta}_{3}\rangle)\!-\! \partial_{x_{3}}S^{\e,\delta}_{2}\right)-\partial^{2}_{x_{3}x_{1}}\!\Psi^{\e,\delta}.
\end{equation}
Now we compute the three functions in \eqref{27a} in terms of the functions $\omed_{i,j}$, we start with $j=1$, 
\begin{align}
\label{27}
&[\partial_{x_{3}}I [\bZed]\times\partial_{x_{3}}\curl\langle\bZed\rangle]_{1}\notag\\
&=\partial_{x_{3}}I [Z^{\e,\delta}_{2}]\partial_{x_{3}}I[\omed_{2,1}]+\partial_{x_{3}}I [Z^{\e,\delta}_{3}]\partial_{x_{3}}I[\omed_{3,1}]\notag\\
&=\partial_{x_{3}}\left(\partial_{x_{3}}I [Z^{\e,\delta}_{2}]I[\omed_{2,1}]\right)-\partial_{x_{3}}^{2}I [Z^{\e,\delta}_{2}]I[\omed_{2,1}]-I[{\dive}_{h}\bZ^{\e,\delta}_{h}]\partial_{x_{3}}I[\omed_{3,1}].
\end{align}

Now, if in the relation \eqref{27} we use  \eqref{22b} and \eqref{23b} and we take into account \eqref{16} and that ${\bf S}^{\e,\delta}$ is bounded in $L^{2}_{t}H^{k}$ for any fixed $k$ and $\delta$, by performing the same type of computation for any $j=1,2,3$ we have, as $\e\to 0$ (for more details see  \cite{GR06} or \cite{FGVN12}),
\begin{equation*}
\label{28}
\frac{1}{\mathbb{T}^{1}}\int_{0}^{T}\!\!\int_{\Omega}\dive(\bZed\otimes\bZed)\varphi dxdt \longrightarrow\frac{1}{\mathbb{T}^{1}}\int_{0}^{T}\!\!\int_{\R^{2}}\dive(\bu^{\delta}\otimes\bu^{\delta})\varphi_{h} dx_{h}dt,
\end{equation*}
for any fixed $\delta>0$ and for any $\varphi\in C^{\infty}_{0}([0,T]\times\R^{2};\R^{2})$, $\varphi=(\varphi_{h}(t,x_{h}),0)$, ${\dive}_{h}\varphi_{h}=0.$

\subsection{Proof of the Theorem \ref{T1}}
In oder to conclude the proof of the Theorem \ref{T1} we just need to  recall that $\bue\cdot\nabla\bue=\dive(\bue\otimes\bue)-\bue\dive\bue$ and to write 
\begin{equation*}
\label{29}
\bue\otimes\bue=(\bue-\bued)\otimes\bue+\bued\otimes(\bue-\bued)+\bued\otimes\bued.
\end{equation*}
and to recall that 
\begin{equation*}
\label{30}
\|\bued-\bue\|_{L^{2}(K)}\leq c\delta\|\nabla\bue\|_{H^{1}(K)}\quad\text{uniformly for $\e>0$.}
\end{equation*}
This means that in the weak formulation the nonlinear term $\bue\otimes\bue$ can be replaced by $\bued\otimes\bued$ that has been analyzed in detail in the previous section. The final step in the proof of the Theorem \ref{T1} is to use in \eqref{weak} a solenoidal test function $\varphi\in C^{\infty}_{0}([0,T]\times\R^{2};\R^{3})$, of the form $\varphi=(\varphi_{h}(t,x_{h}),0)$ and send $\e\to 0$. The only term that deserves some attention is the rotating one, that we handle in the following way,
\begin{align*}
\frac{1}{\e}\int_{0}^{T}&\!\!\int_{\Omega}({\bf g}\times \bue)\cdot\varphi dxdt=\frac{1}{\e}\int_{0}^{T}\!\!\int_{\Omega}({\bf g}\times (\bZ^{\e}+\nabla\Psi^{\e}))\cdot\varphi dxdt\notag\\
 &=\frac{1}{\e}\int_{0}^{T}\!\!\int_{\Omega}({\bf g}\times \langle \bZ^{\e}\rangle)\cdot\varphi dxdt+\frac{1}{\e}\int_{0}^{T}\!\!\int_{\Omega}({\bf g}\times \{ \bZ^{\e}\})\cdot\varphi dxdt=0.
\end{align*}
where we used the fact that $P({\bf g}\times \langle\bZ^{\e}\rangle)=0$ and that $\{ \bZ^{\e}\}=\partial_{x_{3}} I(x)$, with  $\displaystyle{ \int_{\mathbb{T}^{1}}I(x)dx_{3}=0}$, while $\varphi$ depends only on $x_{h}$.

\section{Case $\beta=1/2$}
\label{beta uguale1/2}
In this section we investigate the case $\beta=1/2$, where the system reads as follows
\begin{equation}
\begin{cases}
\displaystyle{\partial_{t}\bue+\frac{1}{\e}({\bf g}\times \bue)+\frac{1}{\e}\nabla \pe=\mu\Delta \bue-\left(\bue\cdot\nabla\right)\bue -\frac{1}{2}(\dive \bue)\bue}\\ \\
\e\partial_{t}\pe+ \dive \bue=0.
\end{cases}
\label{50}
\end{equation}
As  we will see in the next sections, in this case,  the fast rotation due to the  Coriolis force of order $1/\e$ prevails on the low Mach number regime.

\subsection{Preliminary convergence results}
As before, from the  uniform energy bounds \eqref{5}-\eqref{8},  we have 
\begin{equation}
\label{51}
\bue \rightharpoonup \bu \quad \text{weakly in $L^{2}([0,T];H^{1}(\Omega)) $},
\end{equation}
\begin{equation}
\label{51p}
\pe \rightharpoonup \pi \quad \text{$\ast$-weakly in $L^{\infty}([0,T];L^{2}(\Omega)) $}.
\end{equation}
Letting $\e\to 0$ in $\eqref{50}_{2}$ we obtain
\begin{equation}
\label{52}
\dive \bu =0 \quad \text{a.e. in $(0,T)\times \Omega$}.
\end{equation}
Now, if we send $\e$ to zero in $\eqref{50}_{1}$ we have 
\begin{equation}
{\bf g}\times \bu+\nabla \pi=0,
\label{53}
\end{equation}
and similarly as in Section \ref{first conv} we may infer that $\pi$ is independent from the third variable $x_{3}$ and also $\bu_{h}=(u_{1}, u_{2})$ doesn't depend on the vertical variable $x_{3}$, moreover $\dive \bu_{h}=0$. This fact together with the boundary conditions \eqref{3} and the $L^{2}$ bound for $\bu$ yields the conclusion
\begin{equation}
\label{54}
\bu=(\bu_{h}(t, x_{h}),0), \quad u_{3}=0.
\end{equation}

Also for the case $\beta=1/2$, since we are in the framework of a low Mach number limit with ill prepared initial data, in order to establish the convergence of the nonlinear terms in $\eqref{50}_{1}$ we have  to investigate the acoustic waves behavior. This will be done in the next section.

\subsection{Analysis of the acoustic propagator}

We start by writing  the system \eqref{50} as follows,
\begin{equation*}
\begin{cases}
\e\partial_{t}\pe+ \dive \bue=0\\ \\
\displaystyle{\e\partial_{t}\bue+ ({\bf g}\times \bue+\nabla \pe)=\e(\mu\Delta \bue-\left(\bue\cdot\nabla\right)\bue -\frac{1}{2}(\dive \bue)\bue)}.
\end{cases}
\label{ae50}
\end{equation*}
Then, the acoustic system is given by
\begin{equation}
\begin{cases}
\e\partial_{t}\pe+ \dive \bue=0\\ \\
\e\partial_{t}\bue+ ({\bf g}\times \bue+\nabla \pe)=\e \dive \mathbb{G}^{\e}_{1}+\e \mathbb{G}^{\e}_{2},
\end{cases}
\label{ae51}
\end{equation}
 where by taking into account \eqref{6}-\eqref{8} we have  $\mathbb{G}^{\e}_{1}\in L^{2}_{t}L^{2}_{x}$, $\mathbb{G}^{\e}_{2}\in L^{2}_{t}L^{1}_{x}$. Obviously the system \eqref{ae51} has to be read in its weak formulation, namely for any test function $\varphi\in C^{\infty}_{c}([0,T)\times \overline{\Omega})$ it holds
 \begin{equation}
 \label{55}
 \int_{0}^{T}\!\!\int_{\Omega}\left(\e\pe\partial_{t}\varphi+\bue\cdot\nabla\varphi\right)dxdt=-\varepsilon\int_{\Omega}p^{\e}_{0}\varphi(0,\cdot)dx.
 \end{equation}
and for a test function $\varphi\in C^{\infty}_{c}([0,T)\times \overline{\Omega};\R^{3})$  we get
\begin{equation}
 \label{56}
 \begin{split}
 \int_{0}^{T}\!\!\int_{\Omega}\big(\e\bue\partial_{t}\varphi-({\bf g}&\times \bue)\cdot\varphi+\pe\dive\varphi\big)dxdt=\\-\varepsilon\int_{0}^{T}\langle\mathbb{G}^{\e},\varphi\rangle&-\varepsilon\int_{\Omega}\bu^{\e}_{0}\varphi(0,\cdot)dx,
 \end{split}
 \end{equation}
 where
 $$\langle\mathbb{G}^{\e},\varphi\rangle=\int_{\Omega}\left(\mathbb{G}^{\e}_{1}:\nabla\varphi+\mathbb{G}^{\e}_{2}\cdot\varphi\right)dx.$$
To study the behavior of the acoustic system we formally define on $L^{2}(\Omega)\times L^{2}(\Omega;\R^{3})$ the operator $\mathcal{W}$ as 
 \begin{equation}
 \label{57}
 \mathcal{W}\begin{pmatrix}
 p\\
 \bu
 \end{pmatrix}=\begin{pmatrix}
 \dive\bu\\
{\bf g}\times \bu+\nabla p
 \end{pmatrix}.
  \end{equation}

The operator ${\mathcal W}$ is called the {\em acoustic propagator} and the goal of this section is to prove that the  component of the vector $(\pe, \bue)$ orthogonal to the null space of ${\mathcal W}$ decays to zero as $\e\to 0$. Indeed, if this is the case,   the  velocity field  is not affected by the fast oscillations of the acoustic waves since they are killed by the fast decay to zero. In order to achieve this goal we perform a spectral analysis of ${\mathcal W}$.

It is a straightforward computation to deduce that the null space of ${\mathcal W}$  is given by the set
\begin{equation}
\label{61}
\begin{split}
Ker({\mathcal W})=
\bigg\{\!(p, \bu)\mid p=p(x_{h}),& \ \bu=\bu(x_{h}),\\
&{\dive}_{h}\bu_{h}=0,\ \nabla_{h}p=(u_{2},-u_{1})\bigg\}.
\end{split}
\end{equation}
To study the point spectrum  of the operator $\mathcal{W}$ it is more convenient to work in the frequency space. For a function $w$, the Fourier transform $\hat{w}$ with respect to the space variables is defined as,
 $$\hat{w}=\hat{w}(\xi_{h}, k), \quad \xi_{h}=(\xi_{1}, \xi_{2})\in \R^{2},\  k\in \mathbb{Z},$$
 where
 $$\hat{w}(\xi_{h}, k)=\int_{0}^{1}\!\!\int_{\R^{2}}e^{-{\rm i}(\xi_{h}\cdot x_{h}+k\cdot x_{3})}dx_{h}dx_{3}.$$
The  eigenvalues problem for ${\mathcal W}$  is set as follows,
 \begin{equation*}
 \label{58}
 \dive \bu=\lambda p,\qquad {\bf g}\times \bu+\nabla p=\lambda\bu,
 \end{equation*}
 which in Fourier variables has the form
 \begin{equation*}
 \label{59}
 {\rm i}\bigg(\sum_{j=1}^{2}\xi_{j}\hat{\bu}_{j}+k\hat{\bu}_{3}\bigg)-\lambda\hat{p}=0,\qquad {\rm i}(\xi_{1}, \xi_{2}, k)\hat{p}-(\hat{\bu}_{2},-\hat{\bu}_{1},0)-\lambda\hat{\bu}=0.
 \end{equation*}
 After some standard computations we obtain 
 \begin{equation}
 \label{60}
 \lambda^{2}=-\frac{1+|\xi |^{2}+k^{2}\pm\sqrt{(1+|\xi |^{2}+k^{2})^{2}-4k^{2}}}{2}.
 \end{equation}
From \eqref{60} we deduce that the only real eigenvalue is $\lambda=0$ that we obtain for $k=0$ and, as a consequence, we have that the space of eigenvectors of ${\mathcal W}$  coincides with  the $Ker({\mathcal W})$  defined in \eqref{61}.
To prove that the components of $(\pe,\bue)$ orthogonal to $Ker({\mathcal W})$ decay to zero we  use the RAGE theorem that we state in the following form 
(see Cycon et al. \cite[Theorem 5.8]{CyFrKiSi}):

\begin{theorem}
\label{rage}
Let $H$ be a Hilbert space, ${A}: {\cal D}({A}) \subset H \to H$ a
self-adjoint operator, $C: H \to H$ a compact operator, and $P_c$
the orthogonal projection onto  $H_c$,
specifically,
\[
H = H_c \oplus {\rm cl}_H \Big\{ {\rm span} \{ w \in H \ | \ w \ \mbox{an eigenvector of} \ A \} \Big\}.
\]

Then
\begin{equation*}
\label{rage1}
\left\| \frac{1}{\tau} \int_0^\tau \exp(-{\rm i} tA ) C P_c \exp(
{\rm i} tA ) \ dt \right\|_{{\cal L}(H)} \to 0 \ \mbox{as}\ \tau
\to \infty.
\end{equation*}
\end{theorem}

\begin{remark}
If the operator $C$ is non-negative and self-adjoint in $H$, then  we have
\begin{equation}
\frac{1}{T}\int_{0}^{T}\left \langle\exp(-{\rm i} \frac{t}{\e}A )C\exp({\rm i} \frac{t}{\e}A )P_{c}X,Y\right \rangle_{H}dt\leq \omega({\e})\|X\|_{H}\|Y\|_{H},
\label{62}
\end{equation}
where $\omega({\e})\to 0$, as $\e\to 0$.
If we take $Y=P_{c}X$ we get 
\begin{equation}
\frac{1}{T}\int_{0}^{T}\left\|\sqrt{C}\exp({\rm i} \frac{t}{\e}A )P_{c}X,\right\|^{2}_{H}dt\leq \omega({\e})\|X\|_{H}^{2},
\label{63}
\end{equation}
and for any $X\in L^{2}(0,T;H)$ we have
\begin{equation}
\frac{1}{T^{2}}\left\|\sqrt{C}\int_{0}^{t}\exp({\rm i} \frac{t-s}{\e}A )X(s)ds\right\|^{2}_{L^{2}(0,T;H)}dt
\leq \omega({\e})\int_{0}^{T}\|X(s)\|_{H}^{2}ds.
\label{64}
\end{equation}
For more details see \cite{DFN12} or \cite{FGN12}.
\end{remark}
We apply the RAGE Theorem \ref{rage} in the case where  the Hilbert space $H$ is
$$H=H_{M}=\{(p,\bu)\mid\ \hat{p}(\xi_{h}, k)=0,\ \hat{\bu}(\xi_{h}, k)=0\  \text{if} \  |\xi_{h}|+|k|>M\}$$
and the operators $A$, $C$, considered on the space $H_{M}$  are given by
$$A={\rm i}\mathcal{W},\ C[v]=P_{M}[\chi v],\quad \chi\in C^{\infty}_{c}(\Omega),\ 0\leq\chi\leq 1,$$
where 
$$P_{M}:L^{2}(\Omega)\times L^{2}(\Omega;\R^{3}) \longrightarrow H_{M}$$ 
denotes the orthogonal projection into $H_{M}$.
If we denote by $p^{\e}_{M}$, ${\bu}^{\e}_{M}$ the orthogonal projection of $\pe$ and $\bue$ into $H_{M}$ respectively and we apply the projector $P_{M}$ to the acoustic system
\eqref{ae51} we get
\begin{equation}
\e\frac{d}{dt}\begin{pmatrix}
p^{\e}_{M}\\
{\bu}^{\e}_{M}
\end{pmatrix}
+\mathcal{W}\begin{pmatrix}
p^{\e}_{M}\\
{\bu}^{\e}_{M}
\end{pmatrix}=\e\begin{pmatrix}
0\\
\mathbb{G}_{\e,M}
\end{pmatrix},
\label{65}
\end{equation}
and $\mathbb{G}_{\e,M}\in H_{M}$. Moreover, taking into account  \eqref{6} and \eqref{8}, we have the following uniform bound in $\e$
$$\left \|\begin{pmatrix}
0\\
\mathbb{G}_{\e,M}
\end{pmatrix}\right\|_{L^{2}(0,T;H_{M})}\leq c(M).$$
By using Duhamel's formula, the solutions of \eqref{65} are given by,
\begin{equation}
\label{66}
\begin{pmatrix}
p^{\e}_{M}\\
{\bu}^{\e}_{M}
\end{pmatrix}=\exp({\rm i} A\frac{t}{\e} )\begin{pmatrix}
p^{\e}_{M}(0)\\
{\bu}^{\e}_{M}(0)
\end{pmatrix}+\int_{0}^{t}\exp({\rm i} A\frac{t-s}{\e} )\begin{pmatrix}
0\\
\mathbb{G}_{\e,M}
\end{pmatrix}ds
\end{equation}
Now we denote by $Q$ the orthogonal projection into $Ker(\mathcal{W})$,
$$Q:L^{2}(\Omega)\times L^{2}(\Omega;\R^{3}) \longrightarrow Ker(\mathcal{W})$$
and we recall  that the point spectrum of $\mathcal{W}$, hence of the operator $A$ is reduced to $0$. By applying \eqref{63} and \eqref{64} we get
\begin{equation}
\label{67}
Q^{\bot}\begin{pmatrix}
p^{\e}_{M}\\
{\bu}^{\e}_{M}
\end{pmatrix}\to 0\quad \text{in $L^{2}((0,T)\times K;\R^{4})$, as $\e\to 0$,}
\end{equation}
for any compact set $K\subset \overline{\Omega}$ and fixed $M$. 

From \eqref{66} we obtain 
\begin{equation}
\label{68}
Q\begin{pmatrix}
p^{\e}_{M}\\
{\bu}^{\e}_{M}
\end{pmatrix}\to \begin{pmatrix}
\pi_{M}\\
\bu_{M}
\end{pmatrix}\quad \text{in $L^{2}((0,T)\times K;\R^{4})$, as $\e\to 0$,}
\end{equation}
where $\pi$ and $\bu$ are the limits defined in \eqref{51p} and \eqref{51}.

\subsection{Proof of the Theorem \ref{T2}}

Now we are ready to prove the Theorem \ref{T2}. 
By combing together \eqref{67} and \eqref{68} we have
\begin{equation}
\label{69}
P_{M}\bue\rightarrow P_{M}\bu \quad \text{in $L^{2}((0,T)\times K;\R^{3})$,}
\end{equation} 
for any compact set $K\subset \Omega$ and fixed $M$.
Finally \eqref{69} with \eqref{51} and the compact embedding of $W^{1,2}(K)$ in $L^{2}(K)$ gives
\begin{equation}
\label{70}
\bue\rightarrow \bu \quad \text{in $L^{2}((0,T)\times K;\R^{3})$, for any compact set $K\subset \Omega$}.
\end{equation}
Having established the convergences \eqref{70} and \eqref{51p} we can pass into the limit in the weak formulation of \eqref{50}. We take a test function $\psi\in C^{\infty}_{c}([0,T)\times\Omega)$ for the weak formulation of $\eqref{50}_{2}$ and we get
\begin{equation}
\int_{0}^{T}\!\!\int_{\Omega}\left(\pe\partial_{t}\psi +\frac{1}{\e}\bue\cdot\nabla\psi\right )dxdt=
\int_{\Omega}\pe_{0}\psi(0,\cdot)dx,
\label{71}
\end{equation}
For the equation $\eqref{50}_{1}$ we use a test function $\varphi\in  C^{\infty}_{c}([0,T)\times\Omega;\R^{3})$ of the form $\varphi=(\nabla^{\bot}_{h}\psi,0)=(\partial_{x_{2}}\psi, -\partial_{x_{1}}\psi,0)$, $\psi\in C^{\infty}_{c}([0,T)\times\Omega)$, hence we have
\begin{align}
\int_{0}^{T}\!\!\int_{\Omega}\Big(\bue
\partial_{t}\varphi &+\bue\otimes\bue:\nabla\varphi-\frac{1}{2}\bue\dive\bue\cdot\varphi+\frac{1}{\e}({\bf g}\times \bue)\cdot\varphi\Big) dxdt\notag\\
&=\int_{0}^{T}\!\!\int_{\Omega}\nabla\bue\cdot\nabla\varphi dxdt-\int_{\Omega}\bu^{\e}_{0}\cdot \varphi(0,\cdot) dx.
\label{72}
\end{align}
By combining together  \eqref{71} and  \eqref{72} and by performing the limit as $\e\to 0$ we obtain,
\begin{align}
\int_{0}^{T}\!\!\int_{\Omega}\Big(\bu
&\partial_{t}\nabla^{\bot}_{h}\psi+\bu\otimes\bu:\nabla\nabla^{\bot}_{h}\psi+\pi\partial_{t}\psi\Big) dxdt\notag\\
&=\int_{0}^{T}\!\!\int_{\Omega}\nabla\bu\cdot\nabla\nabla^{\bot}_{h} dxdt-\int_{\Omega}\bu_{0}\cdot \nabla^{\bot}_{h}(0,\cdot) +p_{0}\psi(0,\cdot) dx.
\label{73}
\end{align}
Finally, since from \eqref{53}  we have $\bu=\nabla^{\bot}_{h}\pi$  and, recalling that $\bu$ and $\pi$ are independent on the variable $x_{3}$, from \eqref{73} we get
\begin{align*}
\!\!\!\!\int_{0}^{T}\!\!\!\int_{\Omega}&\Big(\nabla^{\bot}_{h}\pi
\partial_{t}\nabla^{\bot}_{h}\psi+\nabla^{\bot}_{h}\pi\otimes\nabla^{\bot}_{h}\pi:\nabla\nabla^{\bot}_{h}\psi+\pi\partial_{t}\psi\Big) dx_{h}dt\notag\\
&=\int_{0}^{T}\!\!\!\int_{\Omega}\nabla\nabla^{\bot}_{h}\pi\cdot\nabla\nabla^{\bot}_{h}\psi dx_{h}dt-\int_{\Omega}\!(\bu_{0}\cdot \nabla^{\bot}_{h}\psi(0,\cdot) +p_{0}\psi(0,\cdot) )dx,
\label{74}
\end{align*}
which is  the equation \eqref{r6} in the sense of distribution.

\end{document}